\documentclass[secthm,seceqn,amsthm,ussrhead,12pt]{amsart}
\usepackage[utf8]{inputenc}
\usepackage[english]{babel}
\usepackage{amssymb,amsmath,amsthm,amsfonts,xcolor,enumerate,hyperref,comment,longtable,cleveref}

\usepackage{times}
\usepackage{cite}
\usepackage{pdflscape}
\usepackage{ulem}
\usepackage[mathcal]{euscript}
\usepackage{tikz}
\usepackage{hyperref}
\usepackage{cancel}
\usepackage{stmaryrd}

\usetikzlibrary{arrows}

\newtheorem{theorem}{Theorem}
\newtheorem{lemma}[theorem]{Lemma}

\newtheorem{remark}[theorem]{Remark}

\newtheorem{definition}[theorem]{Definition}

\newtheorem*{theoremA}{Theorem A}
\newtheorem*{theoremB}{Theorem B}

\newenvironment{Proof}[1][Proof.]{\begin{trivlist}
\item[\hskip \labelsep {\bfseries #1}]}{\flushright
$\Box$\end{trivlist}}

\usepackage{stmaryrd}
\usepackage{xcolor}

\begin{document}
\noindent{\bf
The algebraic and geometric classification of nilpotent right alternative algebras}
\footnote{
The authors thank  the referee of the paper for   constructive comments.
The work was supported by Nazarbayev University Faculty Development Competitive Research Grants N090118FD5341; Nazarbayev University Faculty Development Competitive Research Grants N090118FD5342;  FAPESP 2019/03655-4; CNPq 404649/2018-1; RFBR 20-01-00030; AP08052405 of MES RK.
} 

   \

   {\bf
   Nurlan Ismailov$^{a}$,
   Ivan   Kaygorodov$^{b}$ \&
   Manat Mustafa $^{c}$}

\

{\tiny

$^{a}$ Astana IT University, Nur-Sultan, Kazakhstan

 $^{b}$ CMCC, Universidade Federal do ABC, Santo Andr\'e, Brasil

$^{c}$ Department of Mathematics, Nazarbayev University, Nur-Sultan, 010000, Kazakhstan

\smallskip

   E-mail addresses:

\smallskip
Nurlan Ismailov (nurlan.ismail@gmail.com)

Ivan   Kaygorodov (kaygorodov.ivan@gmail.com)

 Manat Mustafa (manat.mustafa@nu.edu.kz)

}

\

\noindent{\bf Abstract}:
{\it We present algebraic and geometric classifications of the  $4$-dimensional  complex nilpotent right alternative  algebras.
Specifically, we find that, up to isomorphism, there are only $9$ non-isomorphic nontrivial nilpotent right alternative algebras. The corresponding geometric variety has dimension $13$ and it is determined by the Zariski closure of $4$ rigid algebras and one one-parametric family of algebras.
} 

\ 

\noindent {\bf Keywords}:
{\it right alternative algebras, nilpotent algebras, algebraic classification, central extension, geometric classification, degeneration.}

 \ 
 
\noindent {\bf MSC2010}:  17D15,   	17A30, 14D06, 14L30.

 \ 
 
\section*{Introduction}

One of the classical problems in the theory of non-associative algebras is to classify (up to isomorphism) the algebras of dimension $n$ from a certain variety defined by some family of polynomial identities. It is typical to focus on small dimensions, and there are two main directions for the classification: algebraic and geometric. Varieties as associative, Jordan, Lie, Leibniz, Novikov, assosymmetric or Zinbiel algebras have been studied from these two approaches 
(\!\cite{ack, ikm19,  gkk18, kkk18, jkk19, degr3, usefi1, degr1, degr2, ha16,   kv16, mazz79, mazz80} and 
\cite{gkk19, BC99, wolf2,wolf1, ikm19, kkk18, GRH,jkk19, GRH2, kv19, ikv17, ikv18, kppv,   kv16, S90,   klp}, respectively).
Using the classification of all $2$-dimensional algebras \cite{kv16}, it is easily checked that all 2-dimensional right alternative algebras are associative.
In the present paper, we give the algebraic  and geometric classification of
$4$-dimensional complex nilpotent right alternative  algebras,
defined by Albert in 1949  \cite{alb1}.

The variety of right alternative  algebras is defined by the following identity:
\[
\begin{array}{rcl} 
(xy)z - x(yz) &=& -(xz)y + x(zy).
\end{array} \]
It admits the associative, alternative and $(-1,1)$- algebras as a subvariety. 
Albert proved that every semisimple right alternative algebra of characteristic not two is alternative \cite{alb2} and 
after that, 
Thedy showed that a finite dimensional right alternative algebra without a nil ideal is alternative \cite{thedy77}.
But, unfortunately, it was proven that the variety of right alternative algebras does not admit the  Wedderburn principal theorem \cite{thedy78}.
In the infinite case we have a different result:
Mikheev constructed an example of an infinite-dimensional simple right alternative algebra that is not alternative over any field \cite{mikh} and 
it  was finally established that a simple right alternative algebra must be either alternative or nil \cite{skos1}.
The study of simple right alternative superalgebras has a very big progress in recent papers of Pchelintsev and Shashkov (see,
for example, \cite{seroleg} and the reference therein).
Right alternative algebras with some additional identities were studied in papers of Isaev, Pchelintsev and others \cite{isaev}.
For example, Isaev gave the negative answer for the Specht problem in the variety of right alternative algebras \cite{isaev}.
Some combinatorial properties of the variety of right alternative algebras were studied by Umirbaev in \cite{ualbay85,ualbay89}.
On the other side, there is an  interest in the study of 
nilpotent, right nilpotent and solvable right alternative algebras (see, for example, \cite{skos3,ser76,ser13}).
So, Pchelintsev proved an analogue of the well-known Zhevlakov theorem: any right alternative Malcev-admissible  nil algebra of bounded index over a field of characteristic zero is solvable \cite{sergey};
Kaygorodov and Popov proved an analogue of Moens theorem:
a finite-dimensional right alternative algebra over a field of characteristic zero  admitting an invertible Leibniz-derivation is right nilpotent \cite{kp16}.

The key step in our method for algebraically classifying right alternative nilpotent algebras is the calculation of central extensions of smaller algebras.  It comes as no surprise that the central extensions of Lie and non-Lie algebras have been exhaustively studied for years. It is interesting both to describe them and to use them to classify different varieties of algebras \cite{hani,hac16,ss78,klp20,zusmanovich}.
Firstly, Skjelbred and Sund devised a method for classifying nilpotent Lie algebras employing central extensions  \cite{ss78}.
Using this method, all the non-Lie central extensions of  all $4$-dimensional Malcev algebras were described  afterwards \cite{hac16}.
Moreover, the method is especially indicated for the classification of nilpotent algebras and it was used to describe
all the $4$-dimensional nilpotent associative algebras \cite{degr1},
all the $4$-dimensional nilpotent Novikov algebras \cite{kkk18},
all the $4$-dimensional nilpotent bicommutative algebras \cite{kpv19},
all the $5$-dimensional nilpotent Jordan algebras \cite{ha16},
all the $5$-dimensional nilpotent restricted Lie algebras \cite{usefi1}, 
all the $6$-dimensional nilpotent Lie algebras \cite{degr3,degr2},
all the $6$-dimensional nilpotent Tortkara  algebras \cite{gkk18}, 
all the $6$-dimensional anticommutative algebras \cite{kkl19},
and some others.

\section{The algebraic classification of nilpotent right alternative algebras}

\subsection{Method of classification of nilpotent algebras}

The objective of this section is to give an analogue of the Skjelbred-Sund method for classifying nilpotent right alternative algebras. As other analogues of this method were carefully explained in, for example, \cite{hac16}, we will give only some important definitions, and refer the interested reader to the previous sources. We will also employ their notations.

Let $({\bf A}, \cdot)$ be a right alternative  algebra over  $\mathbb C$
and ${\bf V}$ a vector space over ${\mathbb C}$. We define the $\mathbb C$-linear space ${\rm Z^{2}}\left(
\bf A,{\bf V} \right) $  as the set of all  bilinear maps $\theta  \colon {\bf A} \times {\bf A} \longrightarrow {{\bf V}}$
such that
\[ \theta(xy,z)-\theta(x,yz)=-\theta(xz,y)+\theta(x,zy). \]
These maps will be called {\it cocycles}. Consider a
linear map $f$ from $\bf A$ to  ${\bf V}$, and set $\delta f\colon {\bf A} \times
{\bf A} \longrightarrow {{\bf V}}$ with $\delta f  (x,y ) =f(xy )$. Then, $\delta f$ is a cocycle, and we define ${\rm B^{2}}\left(
{\bf A},{{\bf V}}\right) =\left\{ \theta =\delta f\ : f\in {\rm Hom}\left( {\bf A},{{\bf V}}\right) \right\} $, a linear subspace of ${\rm Z^{2}}\left( {\bf A},{{\bf V}}\right) $; its elements are called
{\it coboundaries}. The {\it second cohomology space} ${\rm H^{2}}\left( {\bf A},{{\bf V}}\right) $ is defined to be the quotient space ${\rm Z^{2}}
\left( {\bf A},{{\bf V}}\right) \big/{\rm B^{2}}\left( {\bf A},{{\bf V}}\right) $.

\

Let ${\rm Aut}({\bf A}) $ be the automorphism group of the right alternative algebra ${\bf A} $ and let $\phi \in {\rm Aut}({\bf A})$. Every $\theta \in
{\rm {\rm Z^{2}}}\left( {\bf A},{{\bf V}}\right) $ defines $\phi \theta (x,y)
=\theta \left( \phi \left( x\right) ,\phi \left( y\right) \right) $, with $\phi \theta \in {\rm {\rm Z^{2}}}\left( {\bf A},{{\bf V}}\right) $. It is easily checked that ${\rm Aut}({\bf A})$
acts on ${\rm {\rm Z^{2}}}\left( {\bf A},{{\bf V}}\right) $, and that
 ${\rm B^{2}}\left( {\bf A},{{\bf V}}\right) $ is invariant under the action of ${\rm Aut}({\bf A}).$  
 So, we have that ${\rm Aut}({\bf A})$ acts on ${\rm H^{2}}\left( {\bf A},{{\bf V}}\right)$.

\

Let $\bf A$ be a right alternative  algebra of dimension $m<n$ over  $\mathbb C$, ${{\bf V}}$ a $\mathbb C$-vector
space of dimension $n-m$ and $\theta$ a cocycle, and consider the direct sum ${\bf A}_{\theta } = {\bf A}\oplus {{\bf V}}$ with the
bilinear product `` $\left[ -,-\right] _{{\bf A}_{\theta }}$'' defined by $\left[ x+x^{\prime },y+y^{\prime }\right] _{{\bf A}_{\theta }}=
 xy +\theta(x,y) $ for all $x,y\in {\bf A},x^{\prime },y^{\prime }\in {{\bf V}}$.
It is straightforward that ${\bf A_{\theta}}$ is a right alternative algebra if and only if $\theta \in {\rm Z}^2({\bf A}, {{\bf V}})$; it is  called an $(n-m)$-{\it dimensional central extension} of ${\bf A}$ by ${{\bf V}}$.

We also call the
set ${\rm Ann}(\theta)=\left\{ x\in {\bf A}:\theta \left( x, {\bf A} \right)+ \theta \left({\bf A} ,x\right) =0\right\} $
the {\it annihilator} of $\theta $. We recall that the {\it annihilator} of an  algebra ${\bf A}$ is defined as
the ideal ${\rm Ann}(  {\bf A} ) =\left\{ x\in {\bf A}:  x{\bf A}+ {\bf A}x =0\right\}$. Observe
 that
${\rm Ann}\left( {\bf A}_{\theta }\right) =\big({\rm Ann}(\theta) \cap{\rm Ann}({\bf A})\big)
 \oplus {{\bf V}}$.

\

\begin{definition}
Let ${\bf A}$ be an algebra and $I$ be a subspace of $\operatorname{Ann}({\bf A})$. If ${\bf A}$ is a direct sum (as algebras) of ${\bf A}_0$ and $I$
then $I$ is called an {\it annihilator component} of ${\bf A}$.
\end{definition}
\begin{definition}
A central extension of an algebra $\bf A$ without annihilator component is called a {\it non-split central extension}.
\end{definition}

\

The following result is fundamental to the classification method.

\begin{lemma}
Let ${\bf A}$ be an $n$-dimensional right alternative algebra such that $\dim \ {\rm Ann}({\bf A})=m\neq0$. Then there exists, up to isomorphism, a unique $(n-m)$-dimensional right alternative  algebra ${\bf A}'$ and a bilinear map $\theta \in {\rm Z}^2({\bf A}, {{\bf V}})$ with ${\rm Ann}({\bf A})\cap{\rm Ann}(\theta)=0$, where ${\bf V}$ is a vector space of dimension m, such that ${\bf A} \cong {{\bf A}'}_{\theta}$ and
 ${\bf A}/{\rm Ann}({\bf A})\cong {\bf A}'$.
\end{lemma}

For the proof, we refer the reader to~\cite[Lemma 5]{hac16}.


\

Now, we seek a condition on the cocycles to know when two $(n-m)$-central extensions are isomorphic.
Let us fix a basis $e_{1},\ldots ,e_{s}$ of ${{\bf V}}$, and $
\theta \in {\rm Z^{2}}\left( {\bf A},{{\bf V}}\right) $. Then $\theta $ can be uniquely
written as $\theta \left( x,y\right) =
\displaystyle \sum_{i=1}^{s} \theta _{i}\left( x,y\right) e_{i}$, where $\theta _{i}\in
{\rm Z^{2}}\left( {\bf A},\mathbb C\right) $. It holds that $\theta \in
{\rm B^{2}}\left( {\bf A},{{\bf V}}\right) $\ if and only if all 
$\theta _{i}\in {\rm B^{2}}\left( {\bf A},\mathbb C\right) $, for all $i$; and it also holds that ${\rm Ann}(\theta)={\rm Ann}(\theta _{1})\cap{\rm Ann}(\theta _{2})\ldots \cap{\rm Ann}(\theta _{s})$. 
Furthermore, if ${\rm Ann}(\theta)\cap {\rm Ann}\left( {\bf A}\right) =0$, then ${\bf A}_{\theta }$ has an
annihilator component if and only if $\left[ \theta _{1}\right] ,\left[
\theta _{2}\right] ,\ldots ,\left[ \theta _{s}\right] $ are linearly
dependent in ${\rm H^{2}}\left( {\bf A},\mathbb C\right)$ (see \cite[Lemma 13]{hac16}),
where $[\theta_i]$ is the image of $\theta_i$ in ${\rm H}^2(A,{\mathbb C}).$

\;

Recall that, given a finite-dimensional vector space ${{\bf V}}$ over $\mathbb C$, the {\it Grassmannian} ${\rm G}_{k}\left( {{\bf V}}\right) $ is the set of all $k$-dimensional
linear subspaces of $ {{\bf V}}$. Let ${\rm G}_{s}\left( {\rm H^{2}}\left( {\bf A},\mathbb C\right) \right) $ be the Grassmannian of subspaces of dimension $s$ in
${\rm H^{2}}\left( {\bf A},\mathbb C\right) $.
 For ${\rm W}=\left\langle
\left[ \theta _{1}\right] ,\left[ \theta _{2}\right] ,\dots,\left[ \theta _{s}
\right] \right\rangle \in {\rm G}_{s}\left( {\rm H^{2}}\left( {\bf A},\mathbb C
\right) \right) $ and $\phi \in {\rm Aut}({\bf A})$, define $\phi {\rm W}=\left\langle \left[ \phi \theta _{1}\right]
,\left[ \phi \theta _{2}\right] ,\dots,\left[ \phi \theta _{s}\right]
\right\rangle $. It holds that $\phi {\rm W}\in {\rm G}_{s}\left( {\rm H^{2}}\left( {\bf A},\mathbb C \right) \right) $, and this induces an action of ${\rm Aut}({\bf A})$ on ${\rm G}_{s}\left( {\rm H^{2}}\left( {\bf A},\mathbb C\right) \right) $. We denote the orbit of ${\rm W}\in {\rm G}_{s}\left(
{\rm H^{2}}\left( {\bf A},\mathbb C\right) \right) $ under this action  by ${\rm Orb}({\rm W})$. Let
\[
{\rm W}_{1}=\left\langle \left[ \theta _{1}\right] ,\left[ \theta _{2}\right] ,\dots,
\left[ \theta _{s}\right] \right\rangle ,{\rm W}_{2}=\left\langle \left[ \vartheta
_{1}\right] ,\left[ \vartheta _{2}\right] ,\dots,\left[ \vartheta _{s}\right]
\right\rangle \in {\rm G}_{s}\left( {\rm H^{2}}\left( {\bf A},\mathbb C\right)
\right).
\]
Similarly to~\cite[Lemma 15]{hac16}, in case ${\rm W}_{1}={\rm W}_{2}$, it holds that \[ \bigcap\limits_{i=1}^{s}{\rm Ann}(\theta _{i})\cap {\rm Ann}\left( {\bf A}\right) = \bigcap\limits_{i=1}^{s}
{\rm Ann}(\vartheta _{i})\cap{\rm Ann}( {\bf A}) ,\] 
and therefore the set
\[
{\rm T}_{s}({\bf A}) =\left\{ {\rm W}=\left\langle \left[ \theta _{1}\right] ,
\left[ \theta _{2}\right] ,\dots,\left[ \theta _{s}\right] \right\rangle \in
{\rm G}_{s}\left( {\rm H^{2}}\left( {\bf A},\mathbb C\right) \right) : \bigcap\limits_{i=1}^{s}{\rm Ann}(\theta _{i})\cap{\rm Ann}({\bf A}) =0\right\}
\]
is well defined, and it is also stable under the action of ${\rm Aut}({\bf A})$ (see~\cite[Lemma 16]{hac16}).

\

Now, let ${{\bf V}}$ be an $s$-dimensional linear space and let us denote by
${\rm E}\left( {\bf A},{{\bf V}}\right) $ the set of all non-split $s$-dimensional central extensions of ${\bf A}$ by
${{\bf V}}$. We can write
\[
{\rm E}\left( {\bf A},{{\bf V}}\right) =\left\{ {\bf A}_{\theta }:\theta \left( x,y\right) = \sum_{i=1}^{s}\theta _{i}\left( x,y\right) e_{i} \ \ \text{and} \ \ \left\langle \left[ \theta _{1}\right] ,\left[ \theta _{2}\right] ,\dots,
\left[ \theta _{s}\right] \right\rangle \in {\rm T}_{s}({\bf A}) \right\} .
\]

Finally, we are prepared to state our main result, which can be proved as \cite[Lemma 17]{hac16}.

\begin{lemma}
 Let ${\bf A}_{\theta },{\bf A}_{\vartheta }\in {\rm E}\left( {\bf A},{{\bf V}}\right) $. Suppose that $\theta \left( x,y\right) =  \displaystyle \sum_{i=1}^{s}
\theta _{i}\left( x,y\right) e_{i}$ and $\vartheta \left( x,y\right) =
\displaystyle \sum_{i=1}^{s} \vartheta _{i}\left( x,y\right) e_{i}$.
Then the right alternative algebras ${\bf A}_{\theta }$ and ${\bf A}_{\vartheta } $ are isomorphic
if and only if
$${\rm Orb}\left\langle \left[ \theta _{1}\right] ,
\left[ \theta _{2}\right] ,\dots,\left[ \theta _{s}\right] \right\rangle =
{\rm Orb}\left\langle \left[ \vartheta _{1}\right] ,\left[ \vartheta
_{2}\right] ,\dots,\left[ \vartheta _{s}\right] \right\rangle .$$
\end{lemma}

Then, it exists a bijective correspondence between the set of ${\rm Aut}({\bf A})$-orbits on ${\rm T}_{s}\left( {\bf A}\right) $ and the set of
isomorphism classes of ${\rm E}\left( {\bf A},{{\bf V}}\right) $. Consequently we have a
procedure that allows us, given an right alternative  algebra ${\bf A}$ of
dimension $n$, to construct all non-split central extensions of ${\bf A}$
of dimension $s.$
 This procedure is as follows:
\; \;

{\centerline {\textsl{Procedure}}}
 
\begin{enumerate}
\item For a given [nilpotent] right alternative  algebra ${\bf A}$
of dimension $n$, determine ${\rm H^{2}}( {\bf A},\mathbb {C}) $, ${\rm Ann}({\bf A})$ and ${\rm Aut}({\bf A})$.

\item Determine the set of ${\rm Aut}({\bf A})$-orbits on ${\rm T}_{s}({\bf A}) $.

\item For each orbit, construct the right alternative algebra associated with a
representative of it.
\end{enumerate}

\

\subsection{Notations}
Let ${\bf A}$ be an right alternative algebra and fix
a basis $e_{1},e_{2},\dots,e_{n}$. We define the bilinear form
$\Delta _{ij} \colon {\bf A}\times {\bf A}\longrightarrow \mathbb C$
by $\Delta _{ij}\left( e_{l},e_{m}\right) = \delta_{il}\delta_{jm}$.
Then the set $\left\{ \Delta_{ij}:1\leq i, j\leq n\right\} $ is a basis for the linear space of
the bilinear forms on ${\bf A}$, and in particular, every $\theta \in
{\rm Z^{2}}\left( {\bf A},{\bf V}\right) $ can be uniquely written as $
\theta = \displaystyle \sum_{1\leq i,j\leq n} c_{ij}\Delta _{{i}{j}}$, where $
c_{ij}\in \mathbb C$.
Let us fix the following notations:

\begin{longtable}{lll}
$ {\mathcal R}^{i*}_j$& $\mbox{---}$& $j\mbox{th }i\mbox{-dimensional nilpotent ``non-pure'' right alternative algebra (with identity $xyz=0$});$ \\
${\mathcal R}^i_j$& $\mbox{---}$ & $ j\mbox{th }i\mbox{-dimensional nilpotent ``pure'' right alternative algebra (without identity $xyz=0$)};$  \\
$ {\mathfrak{N}}_i$ & $ \mbox{---}$ & $ i\mbox{-dimensional algebra with zero product};$  \\
$ ({\bf A})_{i,j}$ & $ \mbox{---}$ & $ j\mbox{th }i\mbox{-dimensional central extension of }\bf A.$  \\
\end{longtable}

\subsection{The algebraic classification of  $3$-dimensional nilpotent right alternative algebras}
There are no nontrivial $1$-dimensional nilpotent right alternative algebras, and 
there is only one nontrivial $2$-dimensional nilpotent right alternative algebra
(namely, the non-split central extension of the $1$-dimensional algebra with zero product):
$$\begin{array}{ll llll}
{\mathcal R}^{2*}_{1} &:& (\mathfrak{N}_1)_{2,1} &:& e_1 e_1 = e_2.\\
\end{array}$$

From this algebra, we construct the $3$-dimensional nilpotent right alternative algebra ${\mathcal R}^{3*}_{1}={\mathcal R}^{2*}_{1}\oplus{\mathbb C e_3}.$
Also, the reference \cite{cfk18} gives the description of all central extensions of  ${\mathcal R}^{2*}_{1}$ and $\mathfrak{N}_2$.
Choosing the right alternative algebras between them,
we have the classification of all non-split $3$-dimensional nilpotent right alternative algebras:

\begin{longtable}{ll llllllllllll}
$ {\mathcal R}^{3*}_{2}$  &$ :$ & $ (\mathfrak{N}_2)_{3,1}$  &$ :$&$  e_1 e_1 = e_3$   & $ e_2 e_2=e_3$  \\
$ {\mathcal R}^{3*}_{3}$  &$ :$ & $ (\mathfrak{N}_2)_{3,2}$  &$ :$ & $ e_1 e_2=e_3$ & $ e_2 e_1=-e_3$   \\
$ {\mathcal R}^{3*}_{4}(\alpha)$  &$ :$ & $(\mathfrak{N}_2)_{3,3}$  &$:$& $e_1 e_1 = \alpha e_3$  & $e_2 e_1=e_3$  & $e_2 e_2=e_3$  \\
${\mathcal R}^{3}_{1}$ &$ :$ & $ ({\mathcal R}^{2*}_{1})_{3,1}$  &$:$& $e_1 e_1 = e_2$  & $e_1 e_2=e_3$ & $e_2 e_1=e_3$  \\
\end{longtable}

\subsection{$1$-dimensional central extensions of $3$-dimensional  nilpotent right alternative algebras}
\label{centrext}
\subsubsection{The description of second cohomology space of  $3$-dimensional nilpotent right alternative algebras}

\
In the table below, we give the description of the second cohomology space of  $3$-dimensional nilpotent right alternative algebras.

{\tiny

\begin{longtable}{|l|l|l|l|}
\hline
$\bf A$  & ${\rm Z^{2}}\left( {\bf A}\right) $ & ${\rm B}^2({\bf A})$ & ${\rm H}^2({\bf A})$ \\
\hline
\hline
${\mathcal R}^{3*}_{1} $&  
$ \Big\langle \begin{array}{l}\Delta_{11},\Delta_{12}+\Delta_{21}, \Delta_{13}, \Delta_{31},\Delta_{33} \end{array}\Big\rangle$& 
$ \Big\langle \begin{array}{l}\Delta_{11} \end{array}\Big\rangle$&
$\Big\langle \begin{array}{l}[\Delta_{12}]+ [\Delta_{21}],[\Delta_{13}],[\Delta_{31}],[\Delta_{33}] \end{array}\Big\rangle$ \\
\hline

\hline

${\mathcal R}^{3*}_{2}$ & 
$\Big\langle \begin{array}{l} \Delta_{11},\Delta_{12}, \Delta_{21}, \Delta_{22} \end{array}\Big\rangle$& 
$\Big\langle \begin{array}{l} \Delta_{11} +\Delta_{22} \end{array}\Big\rangle$&  
$\Big\langle \begin{array}{l} [\Delta_{12}], [\Delta_{21}], [\Delta_{22}] \end{array}\Big\rangle$ \\

\hline

${\mathcal R}^{3*}_{3}$ & 
$\Big\langle \begin{array}{l}\Delta_{11},\Delta_{12}, \Delta_{13}, \Delta_{21},\Delta_{22},\Delta_{23} \end{array}\Big\rangle$ & 
$\Big\langle \begin{array}{l}\Delta_{12}-\Delta_{21} \end{array}\Big\rangle$ &  
$\Big\langle \begin{array}{l} [\Delta_{11}], [\Delta_{12}], [\Delta_{22}],[\Delta_{13}],[\Delta_{23}] \end{array}\Big\rangle$ \\

\hline
${\mathcal R}^{3*}_{4}(\alpha\neq0)$ & 
$\Big\langle\begin{array}{l} \Delta_{11},\Delta_{12},\Delta_{21},\Delta_{22} \end{array}\Big\rangle$ &
$\Big\langle\begin{array}{l}\alpha\Delta_{11}+\Delta_{21}+\Delta_{22} \end{array}\Big\rangle$ &
$\Big\langle\begin{array}{l} [\Delta_{12}],  [\Delta_{21}], [\Delta_{22}] \end{array}\Big\rangle$ \\

\hline
${\mathcal R}^{3*}_{4}(0)$ &
$\Big\langle \begin{array}{l} \Delta_{11}, \Delta_{12}, \Delta_{21},\Delta_{22},\Delta_{23}+\Delta_{32}, \\ 
\end{array} \Big\rangle$ &
$\Big\langle\begin{array}{l} \Delta_{21}+\Delta_{22} \end{array}\Big\rangle$ &
$\Big\langle \begin{array}{l}  [\Delta_{11}],[\Delta_{12}], [\Delta_{22}],    [\Delta_{23}]+[\Delta_{32}]   \\ 
\end{array} \Big\rangle $
\\
\hline

${\mathcal R}^{3}_{1}$ & 
$\Big\langle \begin{array}{l} \Delta_{11},\Delta_{12}+\Delta_{21}, \Delta_{13}+\Delta_{22}+\Delta_{31} \end{array}\Big\rangle$ & 
$\Big\langle \begin{array}{l} \Delta_{11} , \Delta_{12}+\Delta_{21} \end{array}\Big\rangle$ &  
$ \Big\langle \begin{array}{l} [\Delta_{13}]+[\Delta_{22}]+[\Delta_{31}] \end{array}\Big\rangle$ \\

\hline

\end{longtable}
}

\begin{remark}
From the description of the cocycles of the algebras ${\mathcal R}^{3*}_{2}$ and ${\mathcal R}^{3*}_{4}(\alpha)_{\alpha\neq0}$, it follows that the 1-dimensional central extensions of these algebras are 2-dimensional central extensions of 2-dimensional nilpotent right alternative algebras. Thanks to \cite{cfk18} there are only   split $2$-dimensional central extensions of $2$-dimensional nilpotent right alternative  algebras.
\end{remark}

\subsubsection{Central extensions of ${\mathcal R}^{3*}_{1}$}

Let us use the following notations: 
\[\nabla_1=[\Delta_{12}]+[\Delta_{21}], \nabla_2=[\Delta_{13}], \nabla_3=[\Delta_{31}], \nabla_4= [\Delta_{33}].
\]
The automorphism group of ${\mathcal R}^{3*}_{1}$ consists of invertible matrices of the form

\[\phi=\left(
                             \begin{array}{ccc}
                               x & 0 & 0   \\
                               u & x^2 & w  \\
                               z & 0 & y                               \end{array}\right)
                               .\]

Since
\[
\phi^T
                           \left(\begin{array}{ccc}
                                0 & \alpha_1 & \alpha_2  \\
                                 \alpha_1 & 0 & 0 \\
                                 \alpha_3 & 0 & \alpha_4 \\
                             \end{array}
                           \right)\phi
                           =\left(\begin{array}{ccc}
                                \alpha^* & \alpha_1^* & \alpha^*_2   \\
                                \alpha^*_1 & 0 & 0 \\
                                \alpha^*_3 & 0 & \alpha^*_4 \\
                             \end{array}\right),\]
we have that the action of ${\rm Aut} ({\mathcal R}^{3*}_{1})$ on the subspace
$\langle  \sum\limits_{i=1}^4\alpha_i \nabla_i \rangle$
is given by
$\langle  \sum\limits_{i=1}^4\alpha^*_i \nabla_i \rangle,$ where

\[
\begin{array}{rclcrcl}
\alpha^*_1&=&\alpha_1x^3 & \ & \alpha^*_2&=&\alpha_1wx+\alpha_2xy+\alpha_4yz;\\
\alpha^*_3&=&\alpha_1wx+\alpha_3xy+\alpha_4yz & \ & \alpha^*_4&=&\alpha_4y^2.
\end{array}\]
We are only interested in cocycles with $\alpha_1\neq0$ and $(\alpha_2,\alpha_3,\alpha_4) \neq (0,0,0)$. 
\begin{enumerate}

\item $\alpha_ 4\neq0$, then:  
\begin{enumerate}

\item if $\alpha_2=\alpha_3$, then setting $x=\frac{1}{\sqrt[3]{\alpha_1}}, y=-\frac{1}{\sqrt{\alpha_4}}, z=-\frac{\alpha_3}{\sqrt[3]{\alpha_1}\alpha_4}$ and $w=0$ we have the representative $\langle \nabla_1+\nabla_4 \rangle.$
\item if $\alpha_2\neq\alpha_3$, then setting $x=\frac{(\alpha_2-\alpha_3)^2}{\alpha_1\alpha_4},$          $y=-\frac{(\alpha_2-\alpha_3)^3}{\alpha_1\alpha_4^2},$ 
$z=-\frac{\alpha_2(\alpha_2-\alpha_3)^2}{\alpha_1\alpha_4^2}$ and $w=0,$ then we have the representative $\langle \nabla_1+\nabla_3+\nabla_4 \rangle.$

\end{enumerate}
\item $\alpha_ 4=0$, then $(\alpha_2,\alpha_3)\neq(0,0)$:  
\begin{enumerate}
\item  if  $\alpha_2\neq\alpha_3$,  then setting $x=\frac{1}{\sqrt[3]{\alpha_1}}, y=-\frac{\sqrt[3]{\alpha_1}}{\alpha_2-\alpha_3}$ and $w=\frac{\alpha_2}{\sqrt[3]{\alpha_1}(\alpha_2-\alpha_3)}$, we have the representative $\langle \nabla_1+\nabla_3\rangle.$ 
\item  if  $\alpha_2=\alpha_3$,  then setting $x=1,$ $y=1$ and $w=-\frac{\alpha_2}{\alpha_1}$, we have the representative $\langle\nabla_1\rangle.$ 

\end{enumerate}
\end{enumerate}
Since cocycle space $\langle\nabla_1\rangle$ defines an algebra with annihilator subspace of dimension 2,  we are only interested in the following representatives for  orbits:
 \[ W_1=\langle \nabla_1+\nabla_4 \rangle, W_2=\langle\nabla_1+\nabla_3+\nabla_4 \rangle, W_3=\langle \nabla_1+\nabla_3\rangle.\] 
We claim that they generate distinct orbits under the action of automorphisms of ${\mathcal R}^{3*}_{1}$. We demonstrate distinctness of orbits ${\rm Orb}(W_1)$ and ${\rm Orb}(W_2)$.
For other cases all calculations are similar and we will omit it. 
Suppose $${\rm Orb}(W_1)={\rm Orb}(W_2).$$ Then there exists an automorphism $\phi\in {\rm Aut}({\mathcal R}^{3*}_{1})$ such that $\phi W_1=W_2$. Then 
$$\phi(W_1)=\langle\phi(\nabla_1+\nabla_4)\rangle,$$ 
where
$$\phi(\nabla_1+\nabla_4)=(u^2+x^3)([\Delta_{12}]+[\Delta_{21}])+$$$$(uw+wx+yz)([\Delta_{13}]+[\Delta_{31}])+2ux^2[\Delta_{22}]+(uw+wx^2)([\Delta_{23}]+[\Delta_{32}])+(2w^2+y^2)[\Delta_{33}].$$
Since $\phi W_1=W_2$, for some $\lambda\in\mathbb C$, we obtain 
$$\nabla_1+\nabla_3+\nabla_4=\lambda\phi(\nabla_1+\nabla_4).$$
However, we have a contradiction in system of equations obtained from the equality. Therefore, ${\rm Orb}(W_1)$ and ${\rm Orb}(W_2)$ are distinct. Similarly, one can check that ${\rm Orb}(W_1)\cap{\rm Orb}(W_3)=\emptyset$ and  
${\rm Orb}(W_2)\cap{\rm Orb}(W_3)=\emptyset.$ Hence, $W_1$, $W_2$ and $W_3$ are representatives for distinct orbits.

Consequently, we have the following 4-dimensional right alternative algebras constructed from ${\mathcal R}^{3*}_{1}:$

\begin{longtable}{lllllllllllllll}
${\mathcal R}^4_{2}$ & $e_1e_1=e_2$ & $ e_1e_2=e_4$ & $ e_2e_1=e_4$  & $ e_3e_3=e_4$\\
${\mathcal R}^4_{3}$ & $e_1e_1=e_2$ & $ e_1e_2=e_4$ & $ e_2e_1=e_4$ & $ e_3e_1=e_4$ & $ e_3e_3=e_4$\\
${\mathcal R}^4_{4}$ & $e_1e_1=e_2$ & $e_1e_2=e_2$ & $ e_2e_1=e_4$ & $ e_3e_1=e_4$\\
\end{longtable}

\subsubsection{Central extensions of ${\mathcal R}^{3*}_{3}$}
Let us use the following notations:
\[ \nabla_1= [\Delta_{11}], \nabla_2=[\Delta_{12}], \nabla_3= [\Delta_{22}], \nabla_4= [\Delta_{13}], \nabla_5= [\Delta_{23}].\]
 The automorphism group of ${\mathcal R}^{3*}_{3}$  consists of invertible matrices of the form
\[\phi=\begin{pmatrix} 
x & y & 0 \\
z & w & 0 \\
t & p & xw-yz 
\end{pmatrix}. \]

\medskip
Since
\[\phi^T\begin{pmatrix}
\alpha_1 & \alpha_2 & \alpha_4\\
0 & \alpha_3 & \alpha_5\\
0 & 0 & 0
\end{pmatrix}\phi=
\begin{pmatrix}
\alpha_1^* & \alpha_2^*-\alpha^* & \alpha_4^*\\
\alpha^* & \alpha_3^* & \alpha_5^*\\
0 & 0 & 0
\end{pmatrix},\]
where

\begin{longtable}{lcl}
$\alpha_1^*$&$=$&$\alpha_1x^2+(\alpha_2x+\alpha_3z)z+(\alpha_4x+\alpha_5z)t$\\
$\alpha_2^*$&$=$&$\alpha_1xy+(\alpha_2x+\alpha_3z)w+(\alpha_4x+\alpha_5z)p+\alpha_1xy+(\alpha_2y+\alpha_3w)z+(\alpha_4y+\alpha_5w)t$\\
$\alpha_3^*$&$=$&$\alpha_1y^2+(\alpha_2y+\alpha_3w)w+(\alpha_4y+\alpha_5w)p$\\
$\alpha_4^*$&$=$&$(\alpha_4x+\alpha_5z)(wx-yz)$\\
$\alpha_5^*$&$=$&$(\alpha_4y+\alpha_5w)(wx-yz).$
\end{longtable}
we obtain that the action of $\operatorname{Aut}\left({\mathcal R}^{3*}_{3} \right)$ on a subspace $ \langle \sum\limits_{i=1}^5 \alpha_i \nabla_i  \rangle$  is given by
$ \langle \sum\limits_{i=1}^5 \alpha_i^* \nabla_i  \rangle.$

We only consider cocycles so that $(\alpha_4,\alpha_5) \neq (0,0).$
It is easy to see that we can suppose $\alpha_5\neq 0.$
\begin{enumerate}

\item 
if $\alpha_3 \alpha_4^2 + \alpha_5 (-\alpha_2 \alpha_4 + \alpha_1 \alpha_5)=0,$
by setting $x=0,$ $y=1,$ $z=1,$ $w=-\alpha_4/\alpha_5,$ $t=-\alpha_3/\alpha_5$ and $ 
p=\frac{2\alpha_3 \alpha_4 - \alpha_2 \alpha_5}{\alpha_5^2}$
we have the representative $\langle \nabla_4 \rangle.$

\item 
if $\alpha_3 \alpha_4^2 + \alpha_5 (-\alpha_2 \alpha_4 + \alpha_1 \alpha_5)\neq 0,$
then by setting
$x=0,$ 
$z=1,$ 
$y=-\frac{\alpha_5^3}{ \alpha_3 \alpha_4^2 + \alpha_5 (-\alpha_2 \alpha_4 + \alpha_1 \alpha_5)},$
$w=\frac{\alpha_4 \alpha_5^2}{  \alpha_3 \alpha_4^2 + \alpha_5 (-\alpha_2 \alpha_4 + \alpha_1 \alpha_5)},$
$t=-\frac{\alpha_3}{\alpha_5}$ and 
$ p=-\frac{\alpha_5^3}{ \alpha_3 \alpha_4^2 + \alpha_5 (-\alpha_2 \alpha_4 + \alpha_1 \alpha_5)},$
we have the representative $\langle \nabla_3+\nabla_4 \rangle.$

\end{enumerate}

 Summarizing, we have the following distinct orbits
\[  \langle\nabla_4 \rangle,  \ \langle \nabla_3+\nabla_4 \rangle. \]

Consequently, we have the following 4-dimensional right alternative algebras constructed from ${\mathcal R}^{3*}_{3}:$

\begin{longtable}{lllllllllllllll}
${\mathcal R}^4_{5}$ & $e_1e_2=e_3$ & $ e_1e_3=e_4$ & $ e_2e_1=-e_3$ \\
${\mathcal R}^4_{6}$ & $e_1e_2=e_3$  & $ e_1e_3=e_4$ & $ e_2e_1=-e_3$ & $ e_2e_2=e_4$\\
\end{longtable}

\subsubsection{Central extensions of ${\mathcal R}^{3*}_{4}(0)$}
Let us use the following notations: 
\[\nabla_1=[\Delta_{11}], \nabla_2=[\Delta_{12}], \nabla_3=[\Delta_{22}], \nabla_4=[\Delta_{23}]+[\Delta_{32}].
\]
The automorphism group of ${\mathcal R}^{3*}_{4}(0)$ consists of invertible matrices of the form

\[\phi=\left(
                             \begin{array}{ccc}
                               x & x-y & 0   \\
                               0 & y & 0  \\
                               z & u & xy                               \end{array}\right)
                               .\]

Since
\[
\phi^T
                           \left(\begin{array}{ccc}
                                \alpha_1 & \alpha_2 & 0  \\
                                 0 & \alpha_3 & \alpha_4 \\
                                 0 & \alpha_4 & 0 \\
                             \end{array}
                           \right)\phi
                           =\left(\begin{array}{ccc}
                                \alpha^*_1 & \alpha^*_2 & 0   \\
                                \alpha^* & \alpha^*+\alpha^*_3 & \alpha^*_4 \\
                                0 &  \alpha^*_4  & 0 \\
                             \end{array}\right),\]
we have that the action of ${\rm Aut} ({\mathcal R}^{3*}_{4}(0))$ on the subspace
$\langle  \sum\limits_{i=1}^4\alpha_i \nabla_i \rangle$
is given by
$\langle  \sum\limits_{i=1}^4\alpha^*_i \nabla_i \rangle,$ where

\[
\begin{array}{rcl}
\alpha^*_1&=&\alpha_1x^2\\
\alpha^*_2&=&\alpha_1(x^2-xy)+\alpha_2xy+\alpha_4yz\\
\alpha^*_3&=&(\alpha_1-\alpha_2)(y^2-xy)+\alpha_3y^2+\alpha_4(2uy-yz)\\
\alpha^*_4&=&\alpha_4xy^2.
\end{array}\]
One can see that the element $\alpha_1\nabla_1+\alpha_2\nabla_2+\alpha_3\nabla_3$ defines a central extension of a two dimensional algebra, therefore we only consider cases when $\alpha_4\neq0.$ 
\begin{enumerate}

\item $\alpha_ 1= 0, \alpha_4\neq0.$ Setting $x=\frac{1}{\alpha_4}, y=1, z=-\frac{\alpha_2}{\alpha_4^2}$ and $u=-\frac{2\alpha_2-\alpha_2\alpha_4+\alpha_3\alpha_4}{2\alpha_4^2},$ we have the representative $\langle \nabla_4 \rangle.$
   
\item $\alpha_ 1\neq0, \alpha_4\neq0.$ Setting $x=\frac{1}{\sqrt{\alpha_1}},y=\frac{\sqrt[4]{\alpha_1}}{\sqrt{\alpha_4}},z=\frac{\alpha_1-\alpha_2-\sqrt[4]{\alpha_1}\sqrt{\alpha_4}}{\sqrt{\alpha_1}\alpha_4}$ and $u=\frac{-\sqrt[4]{\alpha_1^3}(\alpha_1-\alpha_2+\alpha_3)+2(\alpha_1-\alpha_2)\sqrt{\alpha_4}-\sqrt[4]{\alpha_1}\alpha_4}{2\sqrt{\alpha_1\alpha_4^3}},$ we have the representative $\langle \nabla_1+\nabla_4 \rangle.$ 
\end{enumerate}
So we obtain the following new  algebras constructed from  ${\mathcal R}^{3*}_{4}(0):$

\begin{longtable}{lllllllllllllll}
${\mathcal R}^4_{7}$ & $e_2e_1=e_3$ & $e_2e_2=e_3$ & $e_2e_3=e_4$ & $e_3e_2=e_4$\\
${\mathcal R}^4_{8}$ & $ e_1e_1=e_4$ &$ e_2e_1=e_3$ & $e_2e_2=e_3$ &  $ e_2e_3=e_4$ & $ e_3e_2=e_4$\\
\end{longtable} 

\subsubsection{Central extensions of ${\mathcal R}^{3}_{1}$}
Since the cohomology space is one-dimensional, we have only one new algebra constructed from  ${\mathcal R}^{3}_{1}:$
\begin{longtable}{lllllllllllllll} 
${\mathcal R}^4_{9}$ & $e_1 e_1 = e_2$ & $ e_1e_2=e_3$ & $ e_2 e_1=e_3$ & $ e_1 e_3=e_4$ & $  e_2 e_2=e_4$ & $ e_3 e_1=e_4$
\end{longtable}

\subsection{The algebraic classification of $4$-dimensional nilpotent right alternative algebras}

We distinguish two main classes of right alternative algebras: the ``pure'' and the ``non-pure'' ones. By the non-pure ones, we mean those satisfying the identities $(xy)z=0$ and $x(yz)=0$; the pure ones are the rest. 

These ``trivial'' algebras can be considered in many varieties of algebras defined by polynomial identities of degree $3$ (associative, Leibniz, Zinbiel, etc.), and they can be expressed as central extensions of suitable algebras with zero product. Those with dimension $4$ are already classified:  the list of the non-anticommutative ones can be found in \cite{demir}, and there is only one  nilpotent and anticommutative.

Regarding the pure $4$-dimensional nilpotent right alternative  algebras, we have the following theorem, whose proof is  based on the classification of $3$-dimensional nilpotent right alternative  algebras and the results of Section \ref{centrext}.

\begin{theoremA}
Let $\mathcal R$ be a nonzero  $4$-dimensional complex  nilpotent ``pure'' right alternative  algebra.
Then, $\mathcal R$ is isomorphic to one of the algebras listed below.
\end{theoremA}

\begin{longtable}{|l|c|llllll|}
 \hline \multicolumn{1}{|c|}{${\mathcal A}$} & \multicolumn{1}{c|}{ $\mathfrak{Der}$ \  ${\mathcal A}$} & \multicolumn{6}{c|}{Multiplication table} \\ \hline 
\endfirsthead

\hline \multicolumn{1}{|c|}{${\mathcal A}$} & \multicolumn{1}{c|}{ $\mathfrak{Der}$ \  ${\mathcal A}$} & \multicolumn{6}{c|}{Multiplication table} \\ \hline 
\endhead

\hline 
\endlastfoot

${\mathcal R}^4_{1}$ &  6  & $e_1 e_1 = e_2 $& $e_1 e_2=e_3$ & $e_2 e_1=e_3 $&&&\\ \hline

${\mathcal R}^4_{2}$ &  5  &
$e_1 e_1 = e_2$ &  $ e_1e_2=e_4$ & $e_2e_1= e_4$ & $e_3e_3=e_4$&& \\\hline

${\mathcal R}^4_{3}$ &  4  &
$e_1 e_1 = e_2$ &  $ e_1e_2=e_4$ & $e_2e_1= e_4$ & $e_3e_1= e_4$ & $e_3e_3=e_4$ &\\ \hline

${\mathcal R}^4_{4}$  &  5  &
$e_1 e_1 = e_2 $&  $ e_1e_2=e_4$ & $e_2e_1=e_4$ & $e_3e_1=e_4$ &&\\ \hline

${\mathcal R}^4_{5}$  &  4  &
 $e_1 e_2 = e_3$& $e_1 e_3=e_4$  &$e_2e_1=-e_3$  &&&\\ \hline
 
${\mathcal R}^4_{6}$  & 4   &
$e_1 e_2 = e_3$  & $e_1 e_3=e_4$ & $e_2e_1=-e_3$ & $e_2 e_2=e_4$&& \\ \hline
    
${\mathcal R}^4_{7} $ & 4  &
$e_2 e_1 = e_3$ & $e_2e_2=e_3$ &   $e_2 e_3=e_4$ &  $e_3 e_2=e_4$&& \\ \hline 

${\mathcal R}^4_{8} $ &  3  & 
$e_1 e_1=e_4$ &  $e_2 e_1 = e_3 $& $e_2e_2=e_3$ &  $e_2 e_3=e_4$ & $e_3e_2=e_4$& \\ \hline
   
${\mathcal R}^4_{9}$  &  4  & 
$e_1 e_1 = e_2$ &  $e_1e_2=e_3$&  $e_1 e_3=e_4$ &   $e_2 e_1=e_3 $ & $e_2e_2=e_4$ & $e_3e_1=e_4$  \\ \hline
  
\end{longtable}

\begin{remark}
Let $\mathcal R$ be a $4$-dimensional nilpotent non-associative right alternative  algebra.
Then $\mathcal R$ is isomorphic to one algebra from the following list
\[{\mathcal R}^4_{5},\ {\mathcal R}^4_{6},\ {\mathcal R}^4_{7},\ {\mathcal R}^4_{8}.\]
\end{remark}

\begin{remark}
Any 4-dimensional nilpotent right alternative algebra is a (-1,1)-algebra. 
\end{remark}
 
\section{The geometric classification of nilpotent right alternative algebras}

\subsection{Definitions and notation} 
Given an $n$-dimensional vector space ${\bf V}$, the set ${\rm Hom}({\bf V} \otimes {\bf V},{\bf V}) \cong {\bf V}^* \otimes {\bf V}^* \otimes {\bf V}$ 
is a vector space of dimension $n^3$. This space inherits the structure of the affine variety $\mathbb{C}^{n^3}.$ 
Indeed, let us fix a basis $e_1,\dots,e_n$ of ${\bf V}$. Then any $\mu\in {\rm Hom}({\bf V} \otimes {\bf V},{\bf V})$ is determined by $n^3$ structure constants $c_{i,j}^k\in\mathbb{C}$ such that
$\mu(e_i\otimes e_j)=\sum_{k=1}^nc_{i,j}^ke_k$. A subset of ${\rm Hom}({\bf V} \otimes {\bf V},{\bf V})$ is {\it Zariski-closed} if it can be defined by a set of polynomial equations in the variables $c_{i,j}^k$ ($1\le i,j,k\le n$).


The general linear group ${\rm GL}({\bf V})$ acts by conjugation on the variety ${\rm Hom}({\bf V} \otimes {\bf V},{\bf V})$ of all algebra structures on ${\bf V}$:
$$ (g * \mu )(x\otimes y) = g\mu(g^{-1}x\otimes g^{-1}y),$$ 
for $x,y\in {\bf V}$, $\mu\in {\rm Hom}({\bf V} \otimes {\bf V},{\bf V})$ and $g\in {\rm GL}({\bf V})$. Clearly, the ${\rm GL}({\bf V})$-orbits correspond to the isomorphism classes of algebra structures on ${\bf V}$. Let $T$ be a set of polynomial identities which is invariant under isomorphism. Then the subset $\mathbb{L}(T)\subset {\rm Hom}({\bf V} \otimes {\bf V},{\bf V})$ of the algebra structures on ${\bf V}$ which satisfy the identities in $T$ is ${\rm GL}({\bf V})$-invariant and Zariski-closed. It follows that $\mathbb{L}(T)$ decomposes into ${\rm GL}({\bf V})$-orbits. The ${\rm GL}({\bf V})$-orbit of $\mu\in\mathbb{L}(T)$ is denoted by $O(\mu)$ and its Zariski closure by $\overline{O(\mu)}$.

Let ${\bf A}$ and ${\bf B}$ be two $n$-dimensional algebras satisfying the identities from $T$ and $\mu,\lambda \in \mathbb{L}(T)$ represent ${\bf A}$ and ${\bf B}$ respectively.
We say that ${\bf A}$ {\it degenerates to} ${\bf B}$ and write ${\bf A}\to {\bf B}$ if $\lambda\in\overline{O(\mu)}$.
Note that in this case we have $\overline{O(\lambda)}\subset\overline{O(\mu)}$. Hence, the definition of a degeneration does not depend on the choice of $\mu$ and $\lambda$. If ${\bf A}\not\cong {\bf B}$, then the assertion ${\bf A}\to {\bf B}$ 
is called a {\it proper degeneration}. We write ${\bf A}\not\to {\bf B}$ if $\lambda\not\in\overline{O(\mu)}$.

Let ${\bf A}$ be represented by $\mu\in\mathbb{L}(T)$. Then  ${\bf A}$ is  {\it rigid} in $\mathbb{L}(T)$ if $O(\mu)$ is an open subset of $\mathbb{L}(T)$.
Recall that a subset of a variety is called {\it irreducible} if it cannot be represented as a union of two non-trivial closed subsets. A maximal irreducible closed subset of a variety is called an {\it irreducible component}.
It is well known that any affine variety can be represented as a finite union of its irreducible components in a unique way.
The algebra ${\bf A}$ is rigid in $\mathbb{L}(T)$ if and only if $\overline{O(\mu)}$ is an irreducible component of $\mathbb{L}(T)$. 

In the present work we use the methods applied to Lie algebras in \cite{BC99,GRH,GRH2,S90}.
First of all, if ${\bf A}\to {\bf B}$ and ${\bf A}\not\cong {\bf B}$, then $\dim \mathfrak{Der}({\bf A})<\dim \mathfrak{Der}({\bf B})$, where $\mathfrak{Der}({\bf A})$ is the Lie algebra of derivations of ${\bf A}$. We will compute the dimensions of the algebras of derivations and will check the assertion ${\bf A}\to {\bf B}$ only for those ${\bf A}$ and ${\bf B}$ such that $\dim \mathfrak{Der}({\bf A})<\dim \mathfrak{Der}({\bf B})$.

To prove  degenerations, we will construct families of matrices parametrized by $t\in\mathbb{C}^*$. Namely, let ${\bf A}$ and ${\bf B}$ be two algebras represented by the structures $\mu$ and $\lambda$ from $\mathbb{L}(T)$ respectively. Let $e_1,\dots, e_n$ be a basis of ${\bf V}$ and $c_{i,j}^k$ ($1\le i,j,k\le n$) be the structure constants of $\lambda$ in this basis. If there exist $a_i^j(t)\in\mathbb{C}$ ($1\le i,j\le n$, $t\in\mathbb{C}^*$) such that $E_i^t=\sum_{j=1}^na_i^j(t)e_j$ ($1\le i\le n$) form a basis of ${\bf V}$ for any $t\in\mathbb{C}^*$, and the structure constants $c_{i,j}^k(t)$ of $\mu$ in the basis $E_1^t,\dots, E_n^t$ satisfy $\lim\limits_{t\to 0}c_{i,j}^k(t)=c_{i,j}^k$, then ${\bf A}\to {\bf B}$. In this case  $E_1^t,\dots, E_n^t$ is called a {\it parametric basis} for ${\bf A}\to {\bf B}$.

To prove  non-degenerations we will use the following lemma (see \cite{GRH}).

\begin{lemma}\label{main}
Let $\mathcal{B}$ be a Borel subgroup of ${\rm GL}({\bf V})$ and $\mathcal{R}\subset \mathbb{L}(T)$ be a $\mathcal{B}$-stable closed subset.
If ${\bf A} \to {\bf B}$ and ${\bf A}$ can be represented by $\mu\in\mathcal{R}$ then there is $\lambda\in \mathcal{R}$ that represents ${\bf B}$.
\end{lemma}
 
In what follows, each time we need to prove some  non-degeneration $\mu\not\to\lambda$, we will define $\mathcal{R}$ by a set of polynomial equations in structure constants $c_{ij}^k$ in such a way that the structure constants of $\mu$ in the basis $e_1,\dots,e_n$ satisfy these equations. We will omit everywhere the verification of the fact that $\mathcal{R}$ is stable under the action of the subgroup of lower triangular matrices and of the fact that $\lambda\not\in\mathcal{R}$ for any choice of basis of ${\bf V}.$ 
To simplify our equations, we will use the notation $A_i=\langle e_i,\dots,e_n\rangle,\ i=1,\ldots,n$ and write simply $A_pA_q\subset A_r$ instead of $c_{ij}^k=0$ ($i\geq p$, $j\geq q$, $k< r$).

\subsection{Rigid $n$-dimensional nilpotent right alternative  algebra}
As follows from \cite{alb1}, every one-generated right alternative algebra is associative;
and from \cite{kv17}, we conclude that a one-generated algebra can not stay in the orbit closure of a non-one-generated algebra (or a family of non-one-generated algebras).
Then, summarizing we have the following lemma.

\begin{lemma}\label{le}
Any $n$-dimensional one-generated nilpotent right alternative algebra is associative and commutative and it is isomorphic to the following algebra
\[{\mathcal R}^{n} : e_{i}e_{j}=e_{i+j}, \quad 1 \leq i, j \leq{n}, \quad  i+j \leq{n}.\]
The algebra ${\mathcal R}^n$ is rigid in the variety of complex $n$-dimensional nilpotent right alternative algebras.
\end{lemma}

\subsection{The geometric classification of $4$-dimensional nilpotent right alternative algebras}
The main result of the present section is the following theorem.

\begin{theoremB}\label{geobl}
The variety of $4$-dimensional complex nilpotent right alternative algebras has $5$
irreducuble components, defined by rigid algebras $\mathcal{R}_{5}^4, \mathcal{R}_{6}^4, \mathcal{R}_{8}^4, \mathcal{R}_{9}^4$
and the family of algebras $\mathfrak{N}_{3}(\alpha).$

\end{theoremB}

\begin{Proof}
Recall that the full description of the degeneration system of $4$-dimensional trivial right alternative algebras was given in \cite{kppv}.
Using the cited result, we have that the variety of $4$-dimensional trivial right alternative algebras has two irreducible components given by the two following
families of algebras:
$$\begin{array}{lllllll}
\mathfrak{N}_2(\alpha)  & e_1e_1 = e_3 &e_1e_2 = e_4  &e_2e_1 = -\alpha e_3 &e_2e_2 = -e_4 \\

\mathfrak{N}_3(\alpha)  & e_1e_1 = e_4 &e_1e_2 = \alpha e_4  &e_2e_1 = -\alpha e_4 &e_2e_2 = e_4  &e_3e_3 = e_4
\end{array}$$

From Lemma \ref{le}, the algebra $\mathcal{R}_{9}^4$ is rigid.
From the considerations about the dimension of derivations in the previous subsection, it follows that  $\mathcal{R}_{8}^4$ is rigid.
Now, considering conditions ${\mathcal R}$ which are satisfied by the algebra
${\mathcal R}_8^4:$
\[ {\mathcal R}=
\left\{
A_1^2 \subseteq A_3,  A_3^2 = 0,   c_{13}^4=c_{31}^4, c_{23}^4=c_{32}^4  
\right\}. \]
It is easy to see that there is no any basis satisfying the condition ${\mathcal R}$
for algebras $\mathfrak{N}_3(\alpha), {\mathcal R}^4_{5}, {\mathcal R}^4_{6}.$
It follows that ${\mathcal R}^4_{8}\not \to  \mathfrak{N}_3(\alpha), {\mathcal R}^4_{5}, {\mathcal R}^4_{6}$ and they give new three irreducible components with the same dimension.

Therefore, if we prove that these  degenerate to the rest of the nilpotent right alternative algebras, the theorem is proved.
All needed degenerations are given below.

{\tiny 
\begin{longtable}{lclll}

${\mathcal R}^{4}_{9} $&$\to$&${\mathcal R}^{4}_{1}$ &
$E_1^t=e_1 $&
$E_2^t= e_2$ \\
&&& $E_3^t=  e_3 $ &
$E_4^t= t^{-1} e_4$\\

${\mathcal R}^{4}_{8} $&$\to$&${\mathcal R}^{4}_{2}$ &
$E_1^t= -(-1 + t)^3 e_1 + (-1 + t)^2 t e_2 - \frac{(t-1)^4}{2 t} e_3 $&
$E_2^t= (-1 + t)^4 t e_3$ \\
&&& $E_3^t= (-1 + t)^2 t^2 e_2 + \frac{1}{2} (-1 + t)^4 e_3 $ &
$E_4^t= (-1 + t)^6 t^2 e_4$\\

${\mathcal R}^{4}_{8} $&$\to$&${\mathcal R}^{4}_{3}$ &
$E_1^t= (-1 + t)^2 e_1 - (-1 + t) t e_2 + \frac{(-1 + t) (1 - t + t^2)}{ 2 t} e_3 $&
$E_2^t= (-1 + t)^2 t e_3 - (-1 + t)^2 t e_4 $\\
&&& $E_3^t= -(-1 + t) t^2 e_2 + \frac{1}{2} (1 - 2 t + t^3) e_3 $&
$E_4^t= -(-1 + t)^3 t^2 e_4$\\

${\mathcal R}^{4}_{8} $&$\to$&${\mathcal R}^{4}_{4}$ &

$E_1^t= \frac{(t-1)^3}{t} e_1 - (t-1)^2 e_2 + \frac{(t-1 )^3}{2 t} e_3$&
$E_2^t= \frac{(t-1)^4}{t} e_3 - \frac{(t-1)^5}{t^2} e_4 $\\
&&&$E_3^t= -(-1 + t)^2 t e_2 + \frac{1}{2} (-1 + t)^3 e_3 $&
$E_4^t= -\frac{(t-1)^6}{t} e_4$\\

${\mathcal R}^{4}_{8} $&$\to$&${\mathcal R}^{4}_{7}$ &
$E_1^t= \frac{1 - t}{t} e_1 +  e_2$& $E_2^t= \frac{1}{t} e_2 $ \\
&&& $E_3^t= \frac{1}{t^2} e_3$&
$E_4^t= \frac{1}{t^3} e_4$\\

${\mathcal R}^{4}_{8} $&$\to$&${\mathfrak N}_{2}(\alpha)$ &
$E_1^t= -\frac{\alpha (\alpha - t)^3}{(\alpha-1) t} e_1 - 
\frac{\alpha(\alpha - t)^2}{\alpha-1} e_2 + 
\frac{\alpha (\alpha - t)^3 (2 \alpha^2 + t - \alpha (1 + t))}{2 (\alpha-1)^2t^2} e_3 $& 
$E_2^t= \frac{\alpha^2 (\alpha - t)^2}{\alpha-1} e_2 - \frac{\alpha^2 (\alpha - t)^3 (\alpha - t + \alpha t)}{2 (\alpha-1)^2 t^2} e_3 $\\
&&& $E_3^t= \frac{\alpha^3 (\alpha - t)^4}{(\alpha-1)^2 t} e_3 - \frac{ \alpha^4 (\alpha - t)^5}{(\alpha-1)^3 t^2} e_4$ &
$E_4^t= \frac{\alpha^4 (\alpha - t)^6}{(\alpha-1)^3 t^2} e_4$\\

\end{longtable}

}

\end{Proof}

\end{document}